\documentclass{amsart}
\usepackage{amsmath,amssymb}
\newcommand{\pr}{{\mathrm{Pr}}}
\newcommand{\hw}{{\widetilde{H}}}
\newcommand{\per}{{\mathrm{per}}}
\newcommand{\am}{{\mathsf{A}}}
\newcommand{\calf}{{\mathcal{F}}}
\newcommand{\calg}{{\mathcal{G}}}
\newcommand{\calh}{{\mathcal{H}}}
\newcommand{\xdif}{{\mathrm{d}}}

\theoremstyle{definition}

\theoremstyle{remark}

\numberwithin{equation}{section}

\begin{document}

\title{Notes on use of generalized entropies in counting}

\author{Alexey E. Rastegin}

\address{Department of Theoretical Physics, Irkutsk State University,
         Irkutsk, Gagarin Bv. 20, 664003, Russia}
\email{rast@api.isu.ru; alexrastegin@mail.ru}
\thanks{The author is grateful to anonymous referees for valuable comments.}


\begin{abstract}
We address an idea of applying generalized entropies in counting
problems. First, we consider some entropic properties that are
essential for such purposes. Using the $\alpha$-entropies of
Tsallis--Havrda--Charv\'{a}t type, we derive several results
connected with Shearer's lemma. In particular, we derive upper
bounds on the maximum possible cardinality of a family of
$k$-subsets, when no pairwise intersections of these subsets may
coincide. Further, we revisit the Minc conjecture. Our approach
leads to a family of one-parameter extensions of
Br\'{e}gman's theorem. A utility of the obtained bounds is explicitly exemplified. 
\end{abstract}

\maketitle

\section{Introduction}\label{sec1}

The concept of entropy is fundamental in both statistical physics
and information theory. It plays a certain role in applying
information-theoretic ideas to combinatorial problems
\cite{mmt12}. Many results of such a kind were reviewed by
Radhakrishnan \cite{radha01} and Galvin \cite{galvin14}. An
entropy approach is often used in studies of colorings of graphs
\cite{galvin15,tetali12}. Applications of the entropy as a
combinatorial tool are typically based on the Shannon entropy and
its conditional form. Meantime, other entropic functions have
found to be useful in various questions \cite{bengtsson}. The
R\'{e}nyi entropy \cite{renyi61} and the
Tsallis--Havrda--Charv\'{a}t (THC) entropy \cite{havrda,tsallis}
are especially important extensions of the Shannon entropy. In
principle, such entropic functions may have combinatorial or
computational applications. For instances, they both have been
used in global thresholding approach to image processing
\cite{shch05}.

The main goal of this study is to address entropy-based approach
to counting problems with use of the Tsallis--Havrda--Charv\'{a}t
entropies. The paper is organized as follows. In Section
\ref{sec2}, we recall properties of the THC entropies and prove a
useful statement. In Section \ref{sec3}, we obtain THC-entropy
versions of some combinatorial results related to the so-called
Shearer lemma. In particular, we consider an upper estimate for
the maximum possible cardinality of a family of $k$-subsets of the
given set, when subsets obey certain restrictions. In Section
\ref{sec4}, we derive one-parameter family of upper bounds on
permanents of square $(0,1)$-matrices. This family is an extension
of the Br\'{e}gman theorem. We describe an example of utility of
the presented extension.

\section{Definitions and properties of the THC $\alpha$-entropies}\label{sec2}

In this section, we briefly recall definitions of the
Tsallis--Havrda--Charv\'{a}t entropies and related conditional
entropies. Required properties of these entropic functionals are
discussed as well. Let discrete random variable $X$ take values on
the finite set $\Omega_{X}$. The non-extensive entropy of strictly
positive degree $\alpha\neq1$ is defined by \cite{tsallis}
\begin{equation}
H_{\alpha}(X):=\frac{1}{1-\alpha}\left(\sum_{{\,}x\in\Omega_{X}} p(x)^{\alpha}
- 1\right)
. \label{tsaent}
\end{equation}
With the factor $\left(2^{1-\alpha}-1\right)^{-1}$ instead of
$(1-\alpha)^{-1}$, this entropic form was considered by Havrda and
Charv\'{a}t \cite{havrda}. In non-extensive statistical mechanics,
the entropy (\ref{tsaent}) is known as the Tsallis entropy. It is
instructive to rewrite (\ref{tsaent}) as
\begin{equation}
H_{\alpha}(X)=-\sum_{x\in\Omega_{X}}p(x)^{\alpha}{\,}\ln_{\alpha}p(x)
=\sum_{x\in\Omega_{X}}p(x){\,}\ln_{\alpha}\!\left(\frac{1}{p(x)}\right)
. \label{tsaln}
\end{equation}
Assuming $\xi>0$, the so-called $\alpha$-logarithm is defined as
\begin{equation}
\ln_{\alpha}(\xi)=
\left\{\begin{array}{ll}
        \frac{\xi^{1-\alpha}-1}{1-\alpha}\>, & \text{ for } \alpha>0\,,\ \alpha\neq1\,, \\
        \ln\xi\,, & \text{ for } \alpha=1\,.
       \end{array}\right.
\label{lnal}
\end{equation}
In the limit $\alpha\to1$, the entropy (\ref{tsaent}) gives the standard Shannon entropy
\begin{equation}
H_{1}(X)=-\sum_{x\in\Omega_{X}}p(x){\,}\ln{p(x)}
\ . \label{shaln}
\end{equation}
For all real $q\in[0,1]$, we write the binary THC entropy
\begin{equation}
h_{\alpha}(q):=-\,q^{\alpha}\ln_{\alpha}(q)-(1-q)^{\alpha}\ln_{\alpha}(1-q)
\ . \label{binen}
\end{equation}
For $q\in(0,1)$, this function is concave and obeys
$h_{\alpha}(q)=h_{\alpha}(1-q)$. The THC entropies 
succeed some natural properties of the Shannon entropy. The
maximal value of (\ref{tsaent}) is equal to
$\ln_{\alpha}|\Omega_{X}|$ and reached with the uniform
distribution. For $\alpha\geq1$, the joint THC entropy of two
random variables obeys \cite{sf06}
\begin{equation}
H_{\alpha}(X,Y)\leq{H}_{\alpha}(X)+H_{\alpha}(Y)
\ . \label{hxyh}
\end{equation}

In applications of information-theoretic methods, the notion of
conditional entropy is widely used \cite{CT91}. Let us put the
particular functional
\begin{equation}
H_{1}(X|y)=-\sum\nolimits_{x}p(x|y){\,}\ln{p}(x|y)
\ , \nonumber
\end{equation}
in which the sum is taken over $x\in\Omega_{X}$. The entropy of
$X$ conditional on knowing $Y$ is defined as \cite{CT91}
\begin{equation}
H_{1}(X|Y):=\sum\nolimits_{y} p(y){\,}H_{1}(X|y)
=-\sum\nolimits_{x}\sum\nolimits_{y} p(x,y){\,}\ln{p}(x|y)
\ , \label{cshen}
\end{equation}
where $p(x|y)=p(x,y)/p(y)$. When the range of
summation is clear from the context, we will omit symbols such as
$\Omega_{X}$ and $\Omega_{Y}$.

In the literature, two kinds of the conditional THC entropy have
been discussed \cite{sf06}. These forms are respectively inspired
by the two expressions given in (\ref{tsaln}). The first form is
defined as \cite{sf06}
\begin{equation}
H_{\alpha}(X|Y):=\sum\nolimits_{y}p(y)^{\alpha}H_{\alpha}(X|y)
\ , \label{hct1}
\end{equation}
where
\begin{equation}
H_{\alpha}(X|y):=\frac{1}{1-\alpha}\left(\sum\nolimits_{x} p(x|y)^{\alpha}-1\right)
=-\sum\nolimits_{x} p(x|y)^{\alpha}{\,}\ln_{\alpha}p(x|y)
\label{hct0}
\end{equation}
and strictly positive $\alpha\neq1$. The conditional entropy
(\ref{hct1}) is, up to a factor, the quantity originally
introduced by Dar\'{o}czy \cite{ZD70}. For any $\alpha>0$, this
conditional entropy obeys the chain rule written as \cite{ZD70}
\begin{equation}
H_{\alpha}(X,Y)=H_{\alpha}(X|Y)+H_{\alpha}(Y)
\ . \label{tecr2}
\end{equation}
Due to nonnegativity of $H_{\alpha}(X|Y)$, for all $\alpha>0$ we
also have
\begin{equation}
H_{\alpha}(X,Y)\geq{H}_{\alpha}(Y)
\ . \nonumber
\end{equation}
The chain rule (\ref{tecr2}) can be extended to more than two
variables. Up to reordering of random variables, this result is
expressed as \cite{sf06}
\begin{equation}
H_{\alpha}(X_{1},X_{2},\ldots,X_{n})
=\sum\nolimits_{j=1}^{n}H_{\alpha}(X_{j}|X_{j-1},\ldots,X_{1})
\ , \label{ehrl}
\end{equation}
where $\alpha>0$. In the case $\alpha=1$, we obtain the chain rule
with the standard conditional entropy (\ref{cshen}). This property
turns to be very essential in entropic approach to counting. The second form of conditional THC entropy is
introduced as \cite{sf06}
\begin{equation}
\hw_{\alpha}(X|Y):=\sum\nolimits_{y}p(y){\,}H_{\alpha}(X|y)
\ . \label{hct2}
\end{equation}
Although the quantity (\ref{hct2}) does not share the chain rule, it has found use in some questions \cite{sf06,rastkyb}. Its definition is based on the formulation, which seems to be more appropriate in the context of dynamical systems and generalized entropy rates \cite{falniowski15,rastit,skz00}. We also have
$\hw_{\alpha}(X|Y)\leq{H}_{\alpha}(X|Y)$ for $\alpha\in(0,1)$ and
$\hw_{\alpha}(X|Y)\geq{H}_{\alpha}(X|Y)$ for
$\alpha\in(1,\infty)$. For $\alpha=1$, the $\alpha$-entropies (\ref{hct1}) and
(\ref{hct2}) both coincide with (\ref{cshen}).

Using entropic approach in counting, several properties of the
conditional entropy are required. One of these properties is the
chain rule. The standard conditional entropy also satisfies
\begin{equation}
H_{1}(X|Y_{1},\ldots,Y_{n-1},Y_{n})
\leq{H}_{1}(X|Y_{1},\ldots,Y_{n-1})
\ . \label{rlem11}
\end{equation}
Thus, conditioning on more can only reduce the conditional
entropy. This relation is sometimes required in counting
\cite{radha01}. Another very useful property of the standard
conditional entropy is formulated as follows. Let $Y\mapsto{f}(Y)$
be some function, whose domain covers the support of random
variable $Y$. Then we have \cite{radha01}
\begin{equation}
H_{1}\bigl(X\big|f(Y)\bigr)\geq{H}_{1}(X|Y)
\ . \label{stfn}
\end{equation}
We shall now establish analogous properties for the conditional
$\alpha$-entropies.

\newtheorem{t21}{Proposition}
\begin{t21}\label{thm21}
Let $X$ and $Y_{1},\ldots,Y_{n}$ be discrete random variables,
where $n\geq1$. For $\alpha\geq1$, the conditional entropy
(\ref{hct1}) satisfies
\begin{equation}
H_{\alpha}(X|Y_{1},\ldots,Y_{n-1},Y_{n})
\leq{H}_{\alpha}(X|Y_{1},\ldots,Y_{n-1})
\ . \label{rlem1}
\end{equation}
For $\alpha>0$, the conditional entropy (\ref{hct2}) satisfies
\begin{equation}
\hw_{\alpha}(X|Y_{1},\ldots,Y_{n-1},Y_{n})
\leq\hw_{\alpha}(X|Y_{1},\ldots,Y_{n-1})
\ . \label{rlem12}
\end{equation}
Let $Y\mapsto{f}(Y)$ be a function of random variable $Y$. For
$\alpha\geq1$, the conditional entropy (\ref{hct1}) satisfies
\begin{equation}
H_{\alpha}\bigl(X\big|f(Y)\bigr)\geq{H}_{\alpha}(X|Y)
\ . \label{qqfn}
\end{equation}
For $\alpha>0$, the conditional entropy (\ref{hct2}) satisfies
\begin{equation}
\hw_{\alpha}\bigl(X\big|f(Y)\bigr)\geq\hw_{\alpha}(X|Y)
\ . \label{qqfn0}
\end{equation}
\end{t21}

\begin{proof}
The results (\ref{rlem1}) and (\ref{rlem12}) were
proved in \cite{sf06,rastqqt} and \cite{rastit}, respectively. Let
us proceed to (\ref{qqfn}) and (\ref{qqfn0}). Since the standard
case is known, we assume $\alpha\neq1$. To each value $u$ of the
function, we assign the subset $\omega_{u}\subseteq\Omega_{Y}$
such that
\begin{equation}
\omega_{u}:=\bigl\{y:{\>}y\in\Omega_{Y},{\>}f(y)=u\bigr\}
{\ .} \nonumber
\end{equation}
Then the probabilities are written as
\begin{equation}
p(x,u)=\sum_{y\in\omega_{u}}p(x,y)
\ , \qquad
p(u)=\sum_{y\in\omega_{u}}p(y)
{\ .} \label{prov}
\end{equation}
The left-hand side of (\ref{qqfn}) is represented as
\begin{equation}
H_{\alpha}\bigl(X\big|f(Y)\bigr)
=\sum_{u\in\Omega_{f(Y)}} p(u)^{\alpha}H_{\alpha}(X|u)
\ . \label{hafb}
\end{equation}
Replacing $p(u)^{\alpha}$ with $p(u)$, we obtain the expression
for $\hw_{\alpha}\bigl(X\big|f(Y)\bigr)$. For strictly positive
$\alpha\neq1$ and $\xi\geq0$, we introduce the function
\begin{equation}
\eta_{\alpha}(\xi)=\frac{\xi^{\alpha}-\xi}{1-\alpha}
{\ .} \nonumber
\end{equation}
In terms of this function, we now write
\begin{equation}
p(u)^{\alpha}H_{\alpha}(X|u)=
\sum_{x\in\Omega_{X}}p(u)^{\alpha}{\,}\eta_{\alpha}\bigl(p(x|u)\bigr)
\>. \nonumber
\end{equation}
As $\eta_{\alpha}^{\prime\prime}(\xi)\leq0$ for the considered
values of $\alpha$, the function $\xi\mapsto\eta_{\alpha}(\xi)$ is
concave. For fixed $x$ and $u$, we put numbers
$\lambda_{y}=p(y)/p(u)$ and $\xi_{y}=p(x,y)/p(y)=p(x|y)$ such that
\begin{equation}
\sum_{y\in\omega_{u}}\lambda_{y}=\frac{p(u)}{p(u)}=1
\ , \qquad
\sum_{y\in\omega_{u}}\lambda_{y}{\,}\xi_{y}=\frac{p(x,u)}{p(u)}=p(x|u)
\ , \nonumber
\end{equation}
according to (\ref{prov}). By Jensen's inequality, we then obtain
\begin{align}
&p(u){\>}\eta_{\alpha}\bigl(p(x|u)\bigr)\geq
p(u)\sum_{y\in\omega_{u}}\lambda_{y}{\,}\eta_{\alpha}(\xi_{y})
=\sum_{y\in\omega_{u}}p(y){\>}\eta_{\alpha}\bigl(p(x|y)\bigr)
\>. \label{jens0}\\
&p(u)^{\alpha}{\>}\eta_{\alpha}\bigl(p(x|u)\bigr)\geq
p(u)^{\alpha}\sum_{y\in\omega_{u}}\lambda_{y}{\,}\eta_{\alpha}(\xi_{y})
=\sum_{y\in\omega_{u}}p(u)^{\alpha-1}p(y){\>}\eta_{\alpha}\bigl(p(x|y)\bigr)
\>. \label{jens}
\end{align}
Summing (\ref{jens0}) with respect to $x\in\Omega_{X}$, for all
the considered values of $\alpha$ one gets
\begin{equation}
p(u){\,}H_{\alpha}(X|u)\geq
\sum_{y\in\omega_{u}}p(y){\,}H_{\alpha}(X|y)
{\ }. \nonumber
\end{equation}
The latter leads to (\ref{qqfn0}) after summing with respect to
$u\in\Omega_{f(Y)}$. For all $y\in\omega_{u}$ and $\alpha>1$, we
have $p(u)\geq{p}(y)$ and
$p(u)^{\alpha-1}p(y)\geq{p}(y)^{\alpha}$. Combining this with
(\ref{jens}) and summing with respect to $x\in\Omega_{X}$, we
obtain
\begin{equation}
p(u)^{\alpha}H_{\alpha}(X|u)\geq
\sum_{y\in\omega_{u}}p(y)^{\alpha}H_{\alpha}(X|y)
{\ }. \nonumber
\end{equation}
Summing this with respect to $u\in\Omega_{f(Y)}$ completes the
proof of (\ref{qqfn}).
\end{proof}

Note that the standard case $\alpha=1$ of (\ref{qqfn}) and
(\ref{qqfn0}) can be proved by repeating the above reasons with
the concave function $\xi\mapsto-\xi\ln\xi$. In the mentioned
ranges of the parameter, the conditional THC entropies (\ref{hct1})
and (\ref{hct2}) enjoy the property with respect to conditioning
on more. The result (\ref{rlem1}) has allowed to derive the
one-parametric extension of entropic Bell inequalities originally
given in \cite{BC88}. Using (\ref{stfn}), one can deduce a
property useful in entropic approach to Bregman's theorem
\cite{radha97,radha01}. We shall now formulate a similar statement
for the $\alpha$-entropies.

\newtheorem{t22}[t21]{Corollary}
\begin{t22}\label{cor22}
Let the support $\Omega_{Y}$ of random variable $Y$ be partitioned
into $m$ mutually disjoint sets $\omega_{j}$ as
\begin{equation}
\Omega_{Y}=\bigcup_{j=1}^{m}\omega_{j}
\ . \nonumber
\end{equation}
Let $\varpi_{j}\subseteq\Omega_{X}$ be defined as
\begin{equation}
\varpi_{j}:=\bigl\{x:{\>}x\in\Omega_{X},{\>}y\in\omega_{j},{\>}p(x|y)\neq0\bigr\}
\>. \nonumber
\end{equation}
If $\varpi_{j}\neq\varpi_{k}$ for all $j\neq{k}$, then
\begin{align}
H_{\alpha}(X|Y)&\leq\sum_{j=1}^{m} \pr[Y\in\omega_{j}]^{\alpha}{\>}\ln_{\alpha}|\varpi_{j}|
&(1\leq\alpha<\infty)\ , \label{qcor}\\
\hw_{\alpha}(X|Y)&\leq\sum_{j=1}^{m} \pr[Y\in\omega_{j}]{\>}\ln_{\alpha}|\varpi_{j}|
&(0<\alpha<\infty)\ . \label{qcor0}
\end{align}
\end{t22}

\begin{proof}
Let us take the function $Y\mapsto{f}_{\omega}(Y)$ such that
$f_{\omega}(y)=\varpi_{j}$ for each $y\in\omega_{j}$. It then follows
from (\ref{qqfn}) and (\ref{qqfn0}) that
\begin{align}
&H_{\alpha}(X|Y)\leq{H}_{\alpha}\bigl(X\big|f_{\omega}(Y)\bigr)=
\sum_{j=1}^{m} \pr[Y\in\omega_{j}]^{\alpha} H_{\alpha}(X|\varpi_{j})
&(1\leq\alpha<\infty)
\ , \label{qcor1}\\
&\hw_{\alpha}(X|Y)\leq\hw_{\alpha}\bigl(X\big|f_{\omega}(Y)\bigr)=
\sum_{j=1}^{m} \pr[Y\in\omega_{j}]{\,}H_{\alpha}(X|\varpi_{j})
&(0<\alpha<\infty)
\ . \label{qcor01}
\end{align}
The quantity $H_{\alpha}(X|\varpi_{j})$ is represented as
the sum
\begin{equation}
H_{\alpha}(X|\varpi_{j})=\sum_{x\in\varpi_{j}}\eta_{\alpha}\bigl(p(x|\omega_{j})\bigr)
\>. \label{qcor2}
\end{equation}
The sum of $p(x|\omega_{j})$ over $x\in\varpi_{j}$ is equal to
$1$, whence the term (\ref{qcor2}) does not exceed
$\ln_{\alpha}|\varpi_{j}|$. Combining this fact with (\ref{qcor1})
and (\ref{qcor01}) completes the proof. 
\end{proof}

Using Corollary \ref{cor22}, we will obtain upper bounds on 
conditional $\alpha$-entropies in some combinatorial problems. To do so, we have to estimate
not only cardinalities $|\varpi_{j}|$, but also probabilities
$\pr[Y\in\omega_{j}]$. From this viewpoint, the inequality
(\ref{qcor0}) seems to be more appropriate.

\section{Shearer's lemma and intersections of $k$-element sets}\label{sec3}

In this section, we will examine some questions connected with the
Shearer lemma \cite{cfgs86}. The properties of the THC
entropies lead to a lot of inequalities with interesting
combinatorial applications. We first note the following.

\newtheorem{t31}[t21]{Proposition}
\begin{t31}\label{thm31}
Let $X=(X_{1},\ldots,X_{n})$ be a random variable taking values in
the set $S=S_{1}\times\cdots\times{S}_{n}$, where each coordinate
$X_{j}$ is a random variable taking values in $S_{j}$. For all
$\alpha\geq1$, we have
\begin{equation}
H_{\alpha}(X)\leq\sum_{j=1}^{n} H_{\alpha}(X_{j})
\ . \label{hxyh31}
\end{equation}
\end{t31}

\begin{proof}
The claim (\ref{hxyh31}) immediately follows by induction from
(\ref{hxyh}).
\end{proof}

The result (\ref{hxyh31}) is a straightforward extension of
proposition 15.7.2 of the book \cite{nasp08}. Hence, we can obtain
several corollaries. The first of them is posed as follows.

\newtheorem{t32}[t21]{Corollary}
\begin{t32}\label{thm32}
Let $\calf$ be a family of subsets of the set $\{1,\ldots,n\}$,
and let $q_{j}$ denote the fraction of members of $\calf$ that
contain $j$. For all $\alpha\geq1$, we have
\begin{equation}
\ln_{\alpha}|\calf|\leq\sum_{j=1}^{n} h_{\alpha}(q_{j})
\ . \label{hbn1}
\end{equation}
\end{t32}

\begin{proof}
To each set $F\in\calf$, we assign its characteristic vector
$v(F)$, which is a binary $n$-tuple. Let $X=(X_{1},\ldots,X_{n})$
be the random variable taking values in $\{0,1\}^{n}$ such that
\begin{equation}
\pr\bigl[X=v(F)\bigr]=|\calf|^{-1}
\qquad \forall{\,}F\in\calf
\ , \label{fvvf}
\end{equation}
whence $H_{\alpha}(X)=\ln_{\alpha}|\calf|$. The random variable
$X_{j}$ takes values in $\{0,1\}$ and is $j$-th value in the
characteristic vector. By definition of $q_{j}$, the entropy of
$X_{j}$ is equal to $h_{\alpha}(q_{j})$. Combining this with
(\ref{hxyh31}) completes the proof. 
\end{proof}

The result (\ref{hbn1}) provides an upper estimate for the maximum
possible cardinality of a family of subsets. It is an
$\alpha$-entropy version of the basic lemma proved in
\cite{kss81}. The authors of \cite{kss81} used tools of
information theory for studying a family of $k$-subsets, which
satisfy some restrictions. We will further apply (\ref{hbn1}) to a
specific family of $k$-element subsets of the set
$\{1,\ldots,n\}$. Suppose that a family
$\calg=\{G_{1},\ldots,G_{m}\}$ of $m$ subsets of the set
$\{1,\ldots,n\}$ obey the implication
\begin{equation}
{``}\{i,j\}\neq\{s,t\}{"}{\,}\Longrightarrow{\,}
{``}G_{i}\cap{G}_{j}\neq{G}_{s}\cap{G}_{t}{"}
\ . \label{prop2}
\end{equation}
That is, no pairwise intersections of the $k$-subsets may
coincide. We aim to estimate cardinality of this family from
above. Let us begin with an auxiliary result.

\newtheorem{lm3}[t21]{Lemma}
\begin{lm3}\label{hllm}
For $\alpha\in[1,3.67]$, the function
$\lambda\mapsto{h}_{\alpha}(\lambda^{2})/\lambda$ is concave for
$\lambda\in\left[0,1/\sqrt{2}\right]$.
\end{lm3}

\begin{proof}
We left out the case $\alpha=1$, for which the
concavity was reported in \cite{kss81} for all $\lambda\in[0,1]$.
During the proof, we will use the following generalization of
Bernoulli's inequality (see, e.g., section 2.4 of the book
\cite{mitr70}). For $-1<x\neq0$, one has
\begin{align}
(1+x)^{r}&>1+rx &(r\notin[0,1])
\ , \label{gbi0}\\
(1+x)^{r}&<1+rx &(0<r<1)
\ . \label{gbi1}
\end{align}
For $\alpha>1$, we can write the expression
\begin{equation}
(\alpha-1){\>}\frac{h_{\alpha}(\lambda^{2})}{\lambda}=
-\lambda^{2\alpha-1}+\frac{1-(1-\lambda^{2})^{\alpha}}{\lambda}
\ . \label{a122}
\end{equation}
The term $-\lambda^{2\alpha-1}$ is concave with respect to
$\lambda$ for all $\alpha>1$. We will show concavity of the second
term in the right-hand side of (\ref{a122}). Let us use the second
derivative test. For arbitrary function $\xi\mapsto{f}(\xi)$, one
has a general expression
\begin{equation}
\frac{\xdif^{2}}{\xdif\lambda^{2}}
{\>}\frac{f(\lambda^{2})}{\lambda}=
\frac{2}{\lambda^{3}}
{\,}\Bigl(2\xi^{2}f^{\prime\prime}(\xi)-\xi{f}^{\prime}(\xi)+f(\xi)\Bigr)
\ , \label{drsc}
\end{equation}
where $\xi=\lambda^{2}$. Substituting
$f_{\alpha}(\xi)=1-(1-\xi)^{\alpha}$ finally gives
\begin{equation}
2\xi^{2}f_{\alpha}^{\prime\prime}(\xi)-\xi{f}_{\alpha}^{\prime}(\xi)+f_{\alpha}(\xi)=
1+(1-\xi)^{\alpha-2}
\bigl(
-1+c_{1}\xi+c_{2}\xi^{2}
\bigr)
{\,}. \label{rdsc}
\end{equation}
The coefficients in (\ref{rdsc}) are calculated as
\begin{equation}
c_{1}=2-\alpha
\ , \qquad
c_{2}=-1+3\alpha-2\alpha^2=\frac{1}{8}-2\left(\alpha-\frac{3}{4}\right)^{2}
\, . \label{twcf}
\end{equation}
We will show that the quantity (\ref{rdsc}) is not positive for
$\alpha\in[1,3.67]$ and $\lambda\in\left[0,1/\sqrt{2}\right]$. Let
us consider separately the cases $\alpha\in[1,2]$ and
$\alpha\in[2,3.67]$.

Taking the intervals $\alpha\in[1,2]$ and $\xi\in[0,1]$, we have
\begin{equation}
(1-\xi)^{2-\alpha}\leq1-(2-\alpha)\xi=1-c_{1}\xi
\qquad (1\leq\alpha\leq2)
\ . \label{alcn}
\end{equation}
This formula is based on (\ref{gbi1}) with $x=-\xi$ and
$r=2-\alpha$. Due to (\ref{alcn}), we rewrite (\ref{rdsc}) in the
form
\begin{equation}
(1-\xi)^{\alpha-2}
\Bigl(
(1-\xi)^{2-\alpha}
-1+c_{1}\xi+c_{2}\xi^{2}
\Bigr)
\leq(1-\xi)^{\alpha-2}c_{2}\xi^{2}\leq0
\ , \label{c2x2}
\end{equation}
where $c_{2}\leq0$ for $\alpha\in[1,2]$ by (\ref{twcf}). Here,
the concavity takes place for all $\lambda\in[0,1]$.

The case $\alpha\geq2$ is more complicated to analysis. Here, we
introduce the positive parameter $\beta=\alpha-2$. The condition
of negativity of (\ref{rdsc}) is then rewritten as
\begin{equation}
(1-\xi)^{\beta}
\bigl(
1+\beta\xi+\gamma\xi^{2}
\bigr)-1=:F_{\beta}(\xi)\geq0
\ , \nonumber
\end{equation}
where $\gamma=-c_{2}=3+5\beta+2\beta^{2}$. This inequality is to
be proved for $\xi=\lambda^{2}\leq1/2$.

We begin with the case $\beta\in[0,1]$. Using the polynomial
$p_{\beta}(\xi)=1+\beta\xi+\gamma\xi^{2}$, we write the derivative
\begin{equation}
\frac{\xdif{F}_{\beta}}{\xdif\xi}=(1-\xi)^{\beta-1}
\Bigl((1-\xi)p_{\beta}^{\prime}(\xi)-\beta\,p_{\beta}(\xi)
\Bigr)
\,. \nonumber
\end{equation}
Doing simple calculations, we easily obtain
\begin{equation}
(1-\xi)p_{\beta}^{\prime}(\xi)-\beta\,p_{\beta}(\xi)=
\xi\Bigl(\bigl(2\gamma-\beta-\beta^{2}\bigr)-\gamma(2+\beta)\xi
\Bigr)
\,. \label{feriF}
\end{equation}
As $\xi\geq0$, the derivative $\xdif{F}_{\beta}/\xdif\xi$ is not
negative, whenever
\begin{equation}
\xi\leq\frac{2\gamma-\beta-\beta^{2}}{\gamma(2+\beta)}
=\frac{6+9\beta+3\beta^{2}}{(3+5\beta+2\beta^{2})(2+\beta)}
\ . \label{maxxi}
\end{equation}
For $\beta\in[0,1]$, the right-hand side of (\ref{maxxi})
monotonically decreases with $\beta$ from $1$ at $\beta=0$ up to
$0.6$ at $\beta=1$. The condition (\ref{maxxi}) is clearly
satisfied for all $\xi\in[0,1/2]$. Here, the function
$F_{\beta}(\xi)$ does not decrease. Combining this with
$F_{\beta}(0)=0$, we finally get $F_{\beta}(\xi)\geq0$ for all
$\beta\in[0,1]$ and $\xi\in[0,1/2]$.

For $\beta\geq1$, we apply $(1-\xi)^{\beta}\geq1-\beta\xi$ due to
(\ref{gbi0}). Thus, the quantity of interest obeys
\begin{equation}
F_{\beta}(\xi)\geq(1-\beta\xi)
\bigl(
1+\beta\xi+\gamma\xi^{2}
\bigr)-1
\geq
\bigl(\gamma-\beta^{2}\bigr)\xi^{2}-\beta\gamma\xi^{3}
\ . \nonumber
\end{equation}
The latter is not negative, whenever
$\gamma-\beta^{2}\geq\beta\gamma\xi$. Due to $\xi\leq1/2$, we can
focus on the inequality
$2\bigl(\gamma-\beta^{2}\bigr)-\beta\gamma\geq0$, or
\begin{equation}
6+7\beta-3\beta^{2}-2\beta^{3}\geq0
\ . \label{becond}
\end{equation}
Inspecting roots of the polynomial, the condition (\ref{becond})
holds for all $\beta\in[0,1.67]$, though we use it only for
$\beta\in[1,1.67]$. The latter completes the proof for
$\alpha\in[3,3.67]$. 
\end{proof}

We have shown concavity of the function
$\lambda\mapsto{h}_{\alpha}(\lambda^{2})/\lambda$ for
$\alpha\in[1,3.67]$ and $\lambda\in\left[0,1/\sqrt{2}\right]$.
When $\alpha\in[1,2]$, the concavity actually holds for
$\lambda\in[0,1]$. We now formulate the desired estimate as
follows.

\newtheorem{tks1}[t21]{Proposition}
\begin{tks1}\label{tkss}
Let $\calg=\{G_{1},\ldots,G_{m}\}$ be a family containing $m$
$k$-subsets of the set $\{1,\ldots,n\}$, and let $\calg$ satisfy
the property (\ref{prop2}). Let $\lambda_{j}$ denote a proportion
of those members of $\calg$ that contain $j$. If the precondition
\begin{equation}
\lambda_{j}\leq\frac{1}{\sqrt{2}}
\label{lj2}
\end{equation}
holds for all $j\in\{1,\ldots,n\}$, then for all $\alpha\in[1,3.67]$,
\begin{equation}
\ln_{\alpha}\!\binom{m}{2}\leq{k}{\>}\frac{h_{\alpha}(\lambda^{2})}{\lambda}
\ , \qquad
\lambda:=\sum_{j=1}^{n}\frac{\lambda_{j}}{k}{\>}\lambda_{j}
\ . \label{lmlm1}
\end{equation}
\end{tks1}

\begin{proof}
Let us consider pairwise intersections of members of
$\calg$. Each $j\in\{1,\ldots,n\}$ will appear in a proportion
\begin{equation}
q_{j}^{*}=\binom{m}{2}^{-1}\binom{\lambda_{j}m}{2}
=\lambda_{j}^{2}-\frac{\lambda_{j}(1-\lambda_{j})}{m-1}
\ . \nonumber
\end{equation}
In the ratio, the denominator gives the number of all pairwise
intersections; the numerator is the number of those pairwise
intersections that contain $j$. The binary entropy (\ref{binen})
is concave for $q\in(0,1)$ and reaches its maximum at the point
$q=1/2$. Hence, it does not decrease on $(0,1/2)$. Then the
precondition (\ref{lj2}) provides
\begin{equation}
h_{\alpha}(q_{j}^{*})\leq{h}_{\alpha}(\lambda_{j}^{2})
\ . \nonumber
\end{equation}
Combining this with Corollary \ref{thm32}, for $\alpha\geq1$ we obtain
\begin{equation}
\ln_{\alpha}\!\binom{m}{2}\leq\sum_{j=1}^{n}h_{\alpha}(\lambda_{j}^{2})
=k\sum_{j=1}^{n}\frac{\lambda_{j}}{k}{\>}
\frac{h_{\alpha}(\lambda_{j}^{2})}{\lambda_{j}}
\ . \label{lmlm}
\end{equation}
We further use $\sum_{j=1}^{n}\lambda_{j}/k=1$ and concavity of
the function $\lambda\mapsto{h}_{\alpha}(\lambda^{2})/\lambda$.
Combining (\ref{lmlm}) with the Jensen inequality completes the
proof.
\end{proof}

We have obtained an implicit upper bound on $m=|\calg|$ in terms
of $k=|G_{j}|$ and the average proportion $\lambda$ of sets
containing a particular element. Our result is a parametric
extension of one of the statements proved in \cite{kss81}. It also
differs in the following two respects. First, the precondition
(\ref{lj2}) is now imposed. On the other hand, the formula
(\ref{lmlm1}) is more explicit in the sense that no unknown
asymptotically small terms appear. For the prescribed value of
$\lambda\in\left[0,1/\sqrt{2}\right]$, we could optimize a bound
with respect to the parameter $\alpha$. The authors of
\cite{kss81} also consider a family of $k$-sets, in which the
intersection of no two is contained in a third. Such estimates are
connected with one of questions raised by Erd\H{o}s.

The statement of Proposition \ref{thm31} allows a certain
extension. In the case of Shannon entropies, extension of such a
kind has been proved by Shearer \cite{cfgs86}. It is often
referred to as the Shearer lemma \cite{tetali12,kahn01,radha01}.
Its generalization in terms of the THC entropies is posed as
follows.

\newtheorem{ts31}[t21]{Proposition}
\begin{ts31}\label{ths31}
Let $X=(X_{1},\ldots,X_{n})$ be a random variable taking values in
the set $S=S_{1}\times\cdots\times{S}_{n}$, where each coordinate
$X_{j}$ is a random variable taking values in $S_{j}$. For a
subset $I$ of $\{1,\ldots,n\}$, let $X(I)$ denote the random
variable $(X_{j})_{j\in{I}}$. Suppose that $\calg$ is a family
of subsets of $\{1,\ldots,n\}$ and each $j\in\{1,\ldots,n\}$
belongs to at least $k$ members of $\calg$. For $\alpha\geq1$, we
then have
\begin{equation}
k{\,}H_{\alpha}(X)\leq\sum_{G\in\calg} H_{\alpha}\bigl(X(G)\bigr)
\ . \label{hxs31}
\end{equation}
\end{ts31}

\begin{proof}
Following \cite{radha01}, we will apply the chain rule. Using
(\ref{ehrl}), for $\alpha\geq1$ we have
\begin{align}
H_{\alpha}(X)&=\sum_{j=1}^{n} H_{\alpha}\bigl(X_{j}\bigm|(X_{k}:{\,}k<j)\bigr)
\ , \nonumber\\
H_{\alpha}\bigl(X(G)\bigr)
&=\sum_{j\in{G}} H_{\alpha}\bigl(X_{j}\bigm|(X_{k}:{\,}k\in{G},{\,}k<j)\bigr)
\nonumber\\
&\geq\sum_{j\in{G}} H_{\alpha}\bigl(X_{j}\bigm|(X_{k}:{\,}k<j)\bigr)
\ . \label{hxms}
\end{align}
The step (\ref{hxms}) follows from (\ref{rlem1}), since any string
$(X_{k}:{\,}k<j)$ contains more elements than
$(X_{k}:{\,}k\in{G},{\,}k<j)$. Summing (\ref{hxms}) with respect
to all $G\in\calg$ gives
\begin{equation}
\sum_{G\in\calg} H_{\alpha}\bigl(X(G)\bigr)\geq
k\sum_{j=1}^{n} H_{\alpha}\bigl(X_{j}\bigm|(X_{k}:{\,}k<j)\bigr)
\ , \label{xhmx}
\end{equation}
because each $j\in\{1,\ldots,n\}$ belongs to at least $k$ members
of $\calg$.
\end{proof}

The statement of Proposition \ref{ths31} is a THC-entropy
extension of the Shearer lemma. A related geometric picture
was described in \cite{radha01}. Interesting geometric
applications are also discussed in \cite{nasp08}. An immediate
consequence of (\ref{hxs31}) is posed as follows.

\newtheorem{cl31}[t21]{Corollary}
\begin{cl31}\label{tcl31}
Let $N$ be a finite set, and let $\calf$ be a family of subsets of
$N$. Let $\calg=\{G_{1},\ldots,G_{m}\}$ be a family of subsets of
$N$ such that each element of $N$ appears in at least $k$ members
of $\calg$. For each $1\leq{j}\leq{m}$, we define
$\calf_{j}:=\{F\cap{G}_{j}:{\>}F\in\calf\}$. For $\alpha\geq1$, we
then have
\begin{equation}
k{\,}\ln_{\alpha}|\calf|\leq\sum_{j=1}^{m} \ln_{\alpha}|\calf_{j}|
\ . \label{klnf}
\end{equation}
\end{cl31}

For $\alpha=1$, the formula (\ref{klnf}) is reduced to a result
originally proved in \cite{cfgs86}. Some applications of the
latter were also described in \cite{cfgs86}. Of course,
applications of such a kind can further be considered on the base
of (\ref{klnf}). In some cases, a family of one-parameter
relations may give a stronger bound. An explicit example of this
situation is the case of upper bounds on permanents of square
$(0,1)$-matrices.

\section{Upper bounds on permanents of $(0,1)$-matrices}\label{sec4}

In this section, we will derive a family of one-parameter upper
bounds on the permanent of a square $(0,1)$-matrix. The well-known
upper bound on permanents has been conjectured by Minc \cite{minc}
and later proved by Br\'{e}gman \cite{bregman}. Br\'{e}gman's
proof is based on the duality theorem of convex programming and
properties of doubly stochastic matrices. A short elementary proof
of this result was given by Schrijver \cite{as78}. Schrijver also
mentioned an upper bound for permanents of arbitrary nonnegative
matrices. A similar proof with randomization is explained in
\cite{nasp08}. Developing an approach with randomization,
Radhakrishnan presented an entropy-based proof \cite{radha97}. Our
aim is to study the question with use of the THC entropies. First,
we recall preliminary facts. Let
$\am=\bigl[\bigl[a(i,j)\bigr]\bigr]$ be a nonnegative
$n\times{n}$-matrix, and let $S_{n}$ denote the set of all
permutations on $\{1,\ldots,n\}$. The permanent of $\am$ is
defined as
\begin{equation}
\per(\am):=\sum_{\sigma\in{S}_{n}}\prod_{i=1}^{n} a\bigl(i,\sigma(i)\bigr)
\ . \label{perdf}
\end{equation}
We further consider matrices with elements $a(i,j)\in\{0,1\}$. By
$S\subseteq{S}_{n}$, we mean the set of permutations $\sigma$ such
that $a\bigl(i,\sigma(i)\bigr)=1$ for all $i\in\{1,\ldots,n\}$. It
is obvious that $\per(\am)=|S|$. It is assumed that the matrix
contain no rows of only zeros, since otherwise its permanent is
certainly zero. We claim the following.

\newtheorem{t41}[t21]{Proposition}
\begin{t41}\label{thm41}
Let $\am$ be a $n\times{n}$ $(0,1)$-matrix with $\per(\am)\neq0$,
and let $r_{i}\neq0$ be a number of ones in $i$-th row
$(i=1,\ldots,n)$. For all $\alpha\geq1$, the permanent of $\am$
obeys the inequality
\begin{equation}
\ln_{\alpha}\bigl(\per(\am)\bigr)
\leq
\sum_{i=1}^{n}\frac{1}{r_{i}}{\,}\sum_{j=1}^{r_{i}}\ln_{\alpha}(j)
\ . \label{vbth}
\end{equation}
\end{t41}

\begin{proof}
Let $\sigma$ be a random permutation chosen uniformly  from $S$.
We then have the value $H_{\alpha}(\sigma)={\ln_{\alpha}}{|S|}$,
which coincides with the left-hand side of (\ref{vbth}). We will
show that, for $\alpha\geq1$, the entropy $H_{\alpha}(\sigma)$
does not exceed the right-hand side of (\ref{vbth}). Let us choose
a random permutation $\tau\in{S}_{n}$ uniformly. Using the chain
rule (\ref{ehrl}), for each permutation $\tau$ we can write
\begin{align}
H_{\alpha}(\sigma)&=H_{\alpha}\Bigl({\sigma}{\bigl(\tau(1)\bigr)}\Bigr)+
H_{\alpha}\Bigl({\sigma}{\bigl(\tau(2)\bigr)}\Bigm|{\sigma}{\bigl(\tau(1)\bigr)}\Bigr)+\cdots
\nonumber\\
&\quad\cdots+H_{\alpha}\Bigl({\sigma}{\bigl(\tau(n)\bigr)}\Bigm|{\sigma}{\bigl(\tau(1)\bigr)},\ldots,{\sigma}{\bigl(\tau(n-1)\bigr)}\Bigr)
\label{hper1}\\
&\leq\hw_{\alpha}\Bigl({\sigma}{\bigl(\tau(1)\bigr)}\Bigr)+
\hw_{\alpha}\Bigl({\sigma}{\bigl(\tau(2)\bigr)}\Bigm|{\sigma}{\bigl(\tau(1)\bigr)}\Bigr)
+\cdots
\nonumber\\
&\quad\cdots+\hw_{\alpha}\Bigl({\sigma}{\bigl(\tau(n)\bigr)}\Bigm|{\sigma}{\bigl(\tau(1)\bigr)},\ldots,{\sigma}{\bigl(\tau(n-1)\bigr)}\Bigr)
{\,}. \label{hper12}
\end{align}
Here, the second inequality holds for $\alpha\geq1$. To the given
permutation $\tau$ and index $i\in\{1,\ldots,n\}$, we assign the
integer $k(\tau,i)\in\{1,\ldots,n\}$ such that
\begin{equation}
k(\tau,i):=\tau^{-1}(i)
\ , \qquad
{\sigma}{\bigl(\tau(k)\bigr)}=\sigma(i)
\ . \nonumber
\end{equation}
Summing (\ref{hper12}) over all $\tau\in{S}_{n}$, we further
obtain
\begin{align}
|S_{n}|{\,}H_{\alpha}(\sigma)&\leq
\sum_{\tau\in{S}_{n}}\left\{
\hw_{\alpha}\Bigl({\sigma}{\bigl(\tau(1)\bigr)}\Bigr)+
\hw_{\alpha}\Bigl({\sigma}{\bigl(\tau(2)\bigr)}\Bigm|{\sigma}{\bigl(\tau(1)\bigr)}\Bigr)
+\cdots
\right.
\nonumber\\
&\left.
\quad\cdots+\hw_{\alpha}\Bigl({\sigma}{\bigl(\tau(n)\bigr)}\Bigm|{\sigma}{\bigl(\tau(1)\bigr)},\ldots,{\sigma}{\bigl(\tau(n-1)\bigr)}\Bigr)
\right\}
\nonumber\\
&=\sum_{\tau\in{S}_{n}}\sum_{i=1}^{n}
\hw_{\alpha}\Bigl(\sigma(i)\Bigm|{\sigma}{\bigl(\tau(1)\bigr)},\ldots,{\sigma}{\bigl(\tau(k-1)\bigr)}\Bigr)
{\,}. \label{coll}
\end{align}
At the last step, we gather the contributions of different
$\sigma(i)$ separately. For the given $\sigma\in{S}$,
$\tau\in{S}_{n}$, and $i\in\{1,\ldots,n\}$, we define
$R_{i}(\sigma,\tau)$ to be the set of those column indices that
differ from
${\sigma}{\bigl(\tau(1)\bigr)},\ldots,{\sigma}{\bigl(\tau(k-1)\bigr)}$
and give $1$'s in $i$-th row \cite{radha97}. By definition of
$r_{i}$, we have $\bigl|R_{i}(\sigma,\tau)\bigr|\leq{r}_{i}$. Using
(\ref{qcor0}), we then rewrite (\ref{coll}) as
\begin{equation}
|S_{n}|{\,}H_{\alpha}(\sigma)\leq\sum_{i=1}^{n}\sum_{\tau\in{S}_{n}}
\sum_{j=1}^{r_{i}}{\,}\underset{\sigma}{\pr}\Bigl[\bigl|R_{i}(\sigma,\tau)\bigr|=j\Bigr]\ln_{\alpha}(j)
\ . \nonumber
\end{equation}
Dividing this relation by $|S_{n}|$ and taking into account the
uniform distribution of $\tau\in{S}_{n}$, we immediately obtain
\begin{equation}
H_{\alpha}(\sigma)\leq\sum_{i=1}^{n}
\sum_{j=1}^{r_{i}}{\,}\underset{\sigma,\tau}{\pr}\Bigl[\bigl|R_{i}(\sigma,\tau)\bigr|=j\Bigr]\ln_{\alpha}(j)
\ . \label{prsita}
\end{equation}
We now recall a principal observation of \cite{radha97} that
\begin{equation}
\underset{\sigma,\tau}{\pr}\Bigl[{\bigl|R_{i}(\sigma,\tau)\bigr|}=j\Bigr]
=\frac{1}{r_{i}}
\ . \nonumber
\end{equation}
Combining this with (\ref{prsita}) completes the proof.
\end{proof}

The statement of Theorem \ref{thm41} leads to an one-parameter
family of upper bounds on permanents. In the limit
$\alpha\to1^{+}$, the relation (\ref{vbth}) leads to the previous
result \cite{nasp08}
\begin{equation}
\per(\am)\leq\prod_{i=1}^{n} (r_{i}!)^{1/r_{i}}
\ . \label{micon}
\end{equation}
It was conjectured in \cite{minc} and then proved in several ways
\cite{bregman,radha97,as78}. This result can naturally be
reformulated as an upper bound on the number of perfect matchings
in a bipartite graph \cite{galvin14,radha01}.

We now consider a significance of the one-parameter bound
(\ref{vbth}). It is instructive to consider a concrete example.
Let $n\times{n}$-matrix $\am$ have elements $a(1,j)=1$ for all
$j=1,\ldots,n$ and $a(i,j)=\delta(i,j)$ for $i=2,\ldots,n$. That
is, our matrix is obtained from the identity $n\times{n}$-matrix
by filling its first row with ones. We then have $\per(\am)=1$. On
the other hand, one gives $r_{1}=n$ and $r_{2}=\cdots=r_{n}=1$.
Let us compare values of the bounds (\ref{vbth}) and
(\ref{micon}). It is easy to apply (\ref{vbth}) in the case
$\alpha=2$, since $\ln_{2}(\xi)=1-1/\xi$ due to (\ref{lnal}). For
$\alpha=2$, the upper bound (\ref{vbth}) gives
\begin{equation}
1-\frac{1}{\per(\am)}\leq\frac{1}{n}{\,}\sum_{j=1}^{n}\left(1-\frac{1}{j}\right)
=1-\frac{\calh_{n}}{n}
\ . \label{harn1}
\end{equation}
By $\calh_{n}$, we denote the $n$-th harmonic number
\cite{conmath}. It is well known that the asymptotic expansion of
this number for large $n$ is written as \cite{conmath}
\begin{equation}
\calh_{n}=\ln{n}+\gamma+O(1/n)
\ , \nonumber
\end{equation}
where $\gamma$ is the Euler--Mascheroni constant. From
(\ref{harn1}), we immediately obtain
\begin{equation}
\per(\am)\leq\frac{n}{\calh_{n}}=\frac{n}{\ln{n}}\left\{1+O\!\left(\frac{1}{\ln{n}}\right)\right\}
. \label{perharm}
\end{equation}
Substituting the same collection of numbers $r_{i}$ into
(\ref{micon}) gives
\begin{equation}
\per(\am)\leq(n!)^{1/n}=\frac{n}{e}\left\{1+O\!\left(\frac{1}{n}\right)\right\}
. \label{perstr}
\end{equation}
At the last step, we used the Stirling approximation. For
sufficiently large $n$, the upper bound (\ref{perharm}) is
significantly stronger than (\ref{perstr}). On the other hand,
both the bounds are very far from the actual value of permanent.
Nevertheless, our example has shown a relevance of the result
(\ref{vbth}) proved for $\alpha\geq1$.

We can further ask for extending bounds with values
$\alpha\in(0,1)$. The corresponding result can be obtained by an
immediate extension of Schrijver's proof \cite{as78}. We have the
following statement.

\newtheorem{t42}[t21]{Proposition}
\begin{t42}\label{thm42}
Let $\am$ be a $n\times{n}$ $(0,1)$-matrix with $\per(\am)\neq0$,
and let $r_{i}\neq0$ be a number of ones in $i$-th row
$(i=1,\ldots,n)$. For $\alpha\in(0,1)$, the permanent of $\am$
obeys the inequality
\begin{equation}
-\ln_{\alpha}\!\left(\frac{1}{\per(\am)}\right)\leq
\sum_{i=1}^{n}\frac{1}{r_{i}}{\,}\sum_{j=1}^{r_{i}}\ln_{\alpha}(j)
\ . \label{vbth1}
\end{equation}
\end{t42}

\begin{proof}
For convenience, we introduce the function
\begin{equation}
g_{\alpha}(\xi):=\frac{\xi^{\alpha}-\xi}{\alpha-1}
=-\xi{\,}\ln_{\alpha}\!\left(\frac{1}{\xi}\right)
\, , \nonumber
\end{equation}
where $\xi>0$ and $\alpha\in(0,1)$. For these values of $\alpha$, the function
$\xi\mapsto{g}_{\alpha}(\xi)$ is convex. Due to the Jensen
inequality, we have
\begin{equation}
g_{\alpha}\!\left(\frac{1}{r}{\,}\sum_{k=1}^{r} \xi_{k}\right)
\leq\frac{1}{r}{\,}\sum_{k=1}^{r} g_{\alpha}(\xi_{k})
\ . \label{jicn}
\end{equation}
We will prove (\ref{vbth1}) by induction on $n$. For $n=1$, the
result is trivial. Suppose that the claim is already proved for
$(n-1)\times(n-1)$-matrices. The permanent of a
$n\times{n}$-matrix can be decomposed as
\begin{equation}
\per(\am)=\sum_{\substack{k=1 \\ a(i,k)=1}}^{n} \per\bigl(\am(i,k)\bigr)
\ . \label{dcik}
\end{equation}
Here, the submatrix $\am(i,k)$ is obtained from $\am$ by
eliminating the $i$-th row and the $k$-th column. Combining
(\ref{jicn}) with (\ref{dcik}) gives
\begin{equation}
\per(\am){\,}(-1){\,}\ln_{\alpha}\!\left(\frac{r_{i}}{\per(\am)}\right)=
r_{i}{\,}g_{\alpha}\!\left(\frac{\per(\am)}{r_{i}}\right)\leq
\sum_{\substack{k=1 \\ a(i,k)=1}}^{n} g_{\alpha}\Bigl\{\per{\bigl(\am(i,k)\bigr)}\Bigr\}
{\>}. \label{jiik}
\end{equation}
From the definition of the $\alpha$-logarithm, we have the
identity
\begin{equation}
\ln_{\alpha}(r\xi)=\ln_{\alpha}(r)+r^{1-\alpha}\ln_{\alpha}(\xi)
\ . \label{laid}
\end{equation}
Summing (\ref{jiik}) with respect to $i\in\{1,\ldots,n\}$, we
therefore obtain
\begin{align}
&\per(\am)\left\{-\sum_{i=1}^{n}\ln_{\alpha}(r_{i})
-n{\,}\ln_{\alpha}\!\left(\frac{1}{\per(\am)}\right)\right\}
\leq\sum_{i=1}^{n}r_{i}{\,}g_{\alpha}\!\left(\frac{\per(\am)}{r_{i}}\right)
\label{ikk0}\\
&\leq\sum_{i=1}^{n}\sum_{\substack{k=1 \\ a(i,k)=1}}^{n}
\per{\bigl(\am(i,k)\bigr)}{\,}(-1){\,}\ln_{\alpha}\!\left(\frac{1}{\per\bigl(\am(i,k)\bigr)}\right)
\label{ijjk0}\\
&=\sum_{\sigma\in{S}}\sum_{i=1}^{n}
(-1){\,}\ln_{\alpha}\!\left(\frac{1}{\per\bigl\{\am\bigl(i,\sigma(i)\bigr)\bigr\}}\right)
. \label{ijjk1}
\end{align}
To prove (\ref{ikk0}), we used (\ref{laid}) and the relation
$n\leq\sum_{i=1}^{n}r_{i}^{1-\alpha}$ satisfied for
$\alpha\in(0,1)$. To justify (\ref{ijjk1}), we note the following
fact. In the double sum (\ref{ijjk1}), the number of terms from
any pair $(i,k)$ equals the number of those $\sigma\in{S}$ for
which $\sigma(i)=k$. The latter number is
$\per{\bigl(\am(i,k)\bigr)}$ for $a(i,k)=1$, and zero otherwise.
We now apply the induction hypothesis to each
$\per{\bigl\{\am{\bigl(i,\sigma(i)\bigr)}\bigr\}}$ in
(\ref{ijjk1}). The left-hand side of (\ref{ikk0}) is no greater
than
\begin{align}
&\sum_{\sigma\in{S}}\sum_{i=1}^{n}
\left\{
\sum_{\substack{\ell\neq{i} \\ a(\ell,\sigma(i))=0}}
\frac{1}{r_{\ell}}{\,}\sum_{j=1}^{r_{\ell}}\ln_{\alpha}(j)
{\,}+\sum_{\substack{\ell\neq{i} \\ a(\ell,\sigma(i))=1}}
\frac{1}{r_{\ell}-1}{\,}\sum_{j=1}^{r_{\ell}-1}\ln_{\alpha}(j)
\right\}
\nonumber\\
&=\sum_{\sigma\in{S}}\sum_{\ell=1}^{n}
\left\{
\sum_{\substack{i\neq\ell \\ a(\ell,\sigma(i))=0}}
\frac{1}{r_{\ell}}{\,}\sum_{j=1}^{r_{\ell}}\ln_{\alpha}(j)
{\,}+\sum_{\substack{i\neq\ell \\ a(\ell,\sigma(i))=1}}
\frac{1}{r_{\ell}-1}{\,}\sum_{j=1}^{r_{\ell}-1}\ln_{\alpha}(j)
\right\}
\label{stps0}\\
&=\sum_{\sigma\in{S}}\sum_{\ell=1}^{n}
\left\{
\frac{n-r_{\ell}}{r_{\ell}}{\,}\sum_{j=1}^{r_{\ell}}\ln_{\alpha}(j)
{\,}+\sum_{j=1}^{r_{\ell}-1}\ln_{\alpha}(j)
\right\}
. \label{stps1}
\end{align}
In the step (\ref{stps0}), we change an order of summation. The
step (\ref{stps1}) is posed as follows. First, the number of $i$
such that $i\neq\ell$ and ${a}{\bigl(\ell,\sigma(i)\bigr)}=0$ is
equal to $(n-r_{\ell})$. Second, the number of $i$ such that
$i\neq\ell$ and ${a}{\bigl(\ell,\sigma(i)\bigr)}=1$ is equal to
$(r_{\ell}-1)$. These observations allow to compute the sums with
respect to $i$ and get (\ref{stps1}). Adding the term
$\per(\am)\sum_{1\leq\ell\leq{n}}\ln_{\alpha}(r_{\ell})$ to both
(\ref{ikk0}) and (\ref{stps1}), we immediately obtain
\begin{equation}
\per(\am){\,}(-n){\,}\ln_{\alpha}\!\left(\frac{1}{\per(\am)}\right)
\leq\per(\am)\sum_{\ell=1}^{n}\frac{n}{r_{\ell}}{\,}\sum_{j=1}^{r_{\ell}}\ln_{\alpha}(j)
\ . \nonumber
\end{equation}
The latter completes the proof.
\end{proof}

In the limit $\alpha\to1^{-}$, the result (\ref{vbth1}) leads to
the previous result (\ref{micon}). In this regard, it is a proper
extension of (\ref{vbth}) to the parameter range $\alpha\in(0,1)$.
Together, the bounds (\ref{vbth}) and (\ref{vbth1}) cover all
the values $\alpha>0$.


\begin{thebibliography}{60}

\bibitem{nasp08}
Alon, N., Spencer, J.H.: The Probabilistic Method. John Wiley {\&} Sons, New York (2008)

\bibitem{bengtsson}
Bengtsson, I, \.{Z}yczkowski, K.: Geometry of Quantum States: An
Introduction to Quantum Entanglement. Cambridge University Press, Cambridge (2006)

\bibitem{BC88}
Braunstein, S.L., Caves, C.M.: Information-theoretic Bell inequalities. Phys. Rev. Lett. {\bf 61}, 662--665 (1988)

\bibitem{bregman}
Br\'{e}gman, L.M.: Some properties of nonnegative matrices and their permanents. Soviet Math. Dokl. {\bf 14}, 945--949 (1973)

\bibitem{cfgs86}
Chung, F.R.K., Frankl, P., Graham, R., Shearer, J.B.: Some intersection theorems for ordered sets and graphs.
J. Comb. Theory, Ser. A {\bf 43}, 23--37 (1986)

\bibitem{CT91}
Cover, T.M., Thomas, J.A.: Elements of Information Theory. John Wiley {\&} Sons, New York (1991)

\bibitem{ZD70}
Dar\'{o}czy, Z.: Generalized information functions. Inform. Control {\bf 16}, 36--51 (1970)

\bibitem{falniowski15}
Falniowski, F.: Generalized conditional entropy -- determinicity of a process and Rokhlin's formula. Open Sys. Information Dyn. {\bf 22}, 1550025 (2015)

\bibitem{sf06}
Furuichi, S.: Information-theoretical properties of Tsallis entropies. J. Math. Phys. {\bf 47}, 023302 (2006)

\bibitem{galvin14}
Galvin, D.: Three tutorial lectures on entropy and counting. arXiv:1406.7872 [math.CO] (2014)

\bibitem{galvin15}
Galvin, D.: Counting colorings of a regular graph. Graphs Combin. {\bf 31}, 629--638 (2015)

\bibitem{tetali12}
Galvin, D., Tetali, P.: On weighted graph homomorphisms. arXiv:1206.3160 [math.CO] (2012)

\bibitem{conmath}
Graham, R.L., Knuth, D.E., Patashnik, O.: Concrete Mathematics. Addison-Wesley, New York (1994)

\bibitem{havrda}
Havrda, J., Charv\'{a}t, F.: Quantification methods of
classification processes: concept of structural $\alpha$-entropy.
Kybernetika {\bf 3}, 30--35 (1967)

\bibitem{kahn01}
Kahn, J.: An entropy approach to the hard-core model on bipartite graphs. Comb. Prob. Comp. {\bf 10}, 219--237 (2001)

\bibitem{kss81}
Kleitman, D.J., Shearer, J., Sturtevant, D.: Intersections of $k$-element sets. Combinatorica {\bf 1}, 381--384 (1981)

\bibitem{mmt12}
Madiman, M., Marcus, A.W., Tetali, P.: Entropy and set cardinality
inequalities for partition-determined functions. Random Struct.
Algorithms {\bf 40}, 1--26 (2012)

\bibitem{minc}
Minc, H.: Upper bounds for permanents of $(0,1)$-matrices. Bull. Amer. Math. Soc. {\bf 69}, 789--791 (1963)

\bibitem{mitr70}
Mitrinovi\'{c}, D.S.: Analytic Inequalities. Springer-Verlag, Berlin (1970)

\bibitem{radha97}
Radhakrishnan, J.: An entropy proof of Br\'{e}gman's theorem. J. Comb. Theory, Ser. A {\bf 77}, 161--164 (1997)

\bibitem{radha01}
Radhakrishnan, J.: Entropy and counting. In J.C. Misra, ed., Computational Mathematics,
Modelling and Algorithms, Narosa Publishers, New Delhi, 146--168 (2003)

\bibitem{rastkyb}
Rastegin, A.E.: Convexity inequalities for estimating generalized
conditional entropies from below. Kybernetika {\bf 48}, 242--253 (2012)

\bibitem{rastqqt}
Rastegin, A.E.: Tests for quantum contextuality in terms of $q$-entropies. Quantum Inf. Comput. {\bf 14}, 0996--1013 (2014)

\bibitem{rastit}
Rastegin, A.E.: Further results on generalized conditional entropies. RAIRO--Theor. Inf. Appl. {\bf 49}, 67--92 (2015)

\bibitem{renyi61}
R\'{e}nyi, A.: On measures of entropy and information. In J. Neyman, ed., Proc.
4th Berkeley Symposium on Mathematical Statistics and Probability, University of California Press, Berkeley, 547--561 (1961)

\bibitem{as78}
Schrijver, A.: A short proof of Minc's conjecrute. J. Comb. Theory, Ser. A {\bf 25}, 80--83 (1978)

\bibitem{shch05}
Shitong, W., Chung, F.L.: Note on the equivalence
relationship between R\'{e}nyi-entropy based and Tsallis-entropy
based image thresholding. Pattern Recognition Lett. {\bf 26},
2309--2312 (2005)

\bibitem{skz00}
S{\l}omczy\'{n}ski, W., Kwapie\'{n}, J., \.{Z}yczkowski, K.: Entropy computing via integration over fractal measures. Chaos {\bf 10}, 180--188 (2000)

\bibitem{tsallis}
Tsallis, C.: Possible generalization of Boltzmann--Gibbs
statistics. J. Stat. Phys. {\bf 52}, 479--487 (1988)


\end{thebibliography}
\end{document}